\documentclass[a4paper]{amsart} 
 
\usepackage{amsmath,amssymb,amscd,verbatim}

\newcommand{\mc}{\mathcal} 
 
\newcommand{\sub}{\subseteq} 
 
\newcommand{\nsub}{\nsubseteq}

\newcommand{\lra}{\Leftrightarrow} 
 
\newcommand{\longra}{\longrightarrow}

\newcommand{\ra}{\Rightarrow}

\newcommand{\sm}{\setminus} 
 
\newcommand{\al}{\alpha} 
 
\newcommand{\La}{\Lambda} 
 
\newcommand{\be}{\beta}

\newcommand{\spec}{\text{Spec}} 
 
\newcommand{\tspec}{t\text{-Spec}} 
 
\newcommand{\tmax}{t\text{-Max}} 
 
\newcommand{\Max}{\text{Max}}

\newcommand{\op}{\operatorname}

\newtheorem{theorem}{Theorem}[section] 
 
\newtheorem{lemma}[theorem]{Lemma} 
 
\newtheorem{proposition}[theorem]{Proposition}

\newtheorem{corollary}[theorem]{Corollary}

\newtheorem{remark}[theorem]{Remark}

\theoremstyle{definition}

\newtheorem{example}[theorem]{Example}

\begin{document} 
 
 \title{$w$-divisorial domains}



\author{Said El Baghdadi} 
 
\address{Department of Mathematics, Facult\'{e} des Sciences et  
Techniques, P.O. Box 523, Beni Mellal, Morocco} 
 
\email{baghdadi@fstbm.ac.ma}

\author{Stefania Gabelli} 
 
\address{Dipartimento di Matematica, Universit\`{a} degli Studi Roma  
Tre, 
Largo S.  L.  Murialdo, 1, 00146 Roma, Italy}

\email{gabelli@mat.uniroma3.it}

%%%%%%%%%%%%%%%%%%%%%%%%%%%%%%%%%%%%%%%%%% 
 
% ABSTRACT 
 
%%%%%%%%%%%%%%%%%%%%%%%%%%%%%%%%%%%%%%%%%%% 

\begin{abstract} We study the class of domains in which each $w$-ideal is 
divisorial, extending several properties of divisorial and totally 
divisorial domains to a much wider class of domains. In particular we 
consider $PvMD$s and Mori domains.

\end{abstract}

\maketitle

%%%%%%%%%%%%%%%%%%%%%%%%%% INTRODUCTION %%%%%%%%%%%%%%%%%%%%%%%%%%%% 

\section*{Introduction} 
 
 The class of domains in which each  
 nonzero ideal is divisorial has been studied, independently and with 
 different methods, by H. Bass  
 \cite{Ba}, E. Matlis \cite{Ma1} and W. Heinzer \cite{H} in the sixties.  
 Following S. Bazzoni 
and L. Salce \cite{BS, B2}, these domains are now called {\it divisorial domains}.  
Among other results, Heinzer proved that an integrally closed domain is 
divisorial if and only if it is a  Pr\"ufer  domain with certain finiteness 
properties \cite[Theorem 5.1]{H}.  
 
Twenty years later E. Houston and M. Zafrullah introduced in \cite {HZ} 
 the class of domains in which each $t$-ideal is divisorial, which  
they called $TV$-{\it domains}, and characterized $PvMD$s with  
this property \cite[Theorem 3.1]{HZ}. 
However they observed that an integrally closed 
$TV$-domain need not be a $PvMD$ \cite[Remark 3.2]{HZ}; thus in some  
sense the class  
of $TV$-domains is not the right setting for extending to $PvMD$s the  
properties of divisorial  Pr\"ufer  domains.  
 
The purpose of this paper is to investigate $w$-{\it divisorial domains}, 
that is domains in which each $w$-ideal is divisorial.  
This class of domains proves to be the most suitable 
 $t$-analogue of divisorial domains. In fact,  
by using this concept we are  
able to improve and generalize several results proved for Noetherian and 
 Pr\"ufer divisorial domains in \cite{BS, H, O, Q}.

The main result of Section 1 is Theorem \ref{Loc}. It states that $R$ is  
a $w$-divisorial domain if and only if $R$ is a weakly Matlis domain  
(that is a  
domain with $t$-finite 
character such that each $t$-prime ideal is contained in a unique  
$t$-maximal ideal) and $R_M$ is a divisorial domain, for each  
$t$-maximal ideal $M$.  In this way we recover the characterization 
of divisorial domains given in \cite[Proposition 5.4]{BS}. 
 
In Section 2, we study the transfer of the properties of  
$w$-divisoriality and divisoriality 
to certain (generalized) rings of fractions, such as localizations at  
($t$-)prime ideals, 
($t$-)flat overrings and ($t$-)subintersections. 
 
In Section 3 we consider $w$-divisorial $PvMD$s.  We prove that $R$ is  
an integrally  
closed  $w$-divisorial domain if and only if 
$R$ is a weakly Matlis $PvMD$ and each $t$-maximal ideal is  
$t$-invertible (Theorem \ref{ic}). This is the $t$-analogue of \cite  
[Theorem 5.1]{H}. 
We also prove that when $R$ is integrally  
closed, each $t$-linked overring of $R$ is $w$-divisorial  
if and only if $R$ is a generalized Krull domain and each $t$-prime  
ideal is contained in a unique $t$-maximal ideal 
(Theorem \ref{twic}). Since in the  Pr\"ufer case generalized Krull domains 
coincide with generalized Dedekind domains  
\cite{elb1}, we obtain that an integrally closed domain is totally 
divisorial  if and only if it is a divisorial generalized Dedekind 
domain \cite[Section 4]{O}. 
 
The last section is devoted to Mori $w$-divisorial domains. A Mori  
$w$-divisorial domain is necessarily of $t$-dimension one and each of its 
localizations at a height-one prime is Noetherian (Corollary \ref{SM}).  
Noetherian divisorial and totally divisorial domains 
were intensely studied in \cite {BS, Ba, Ma1, Q}.  
It turns out that several of the results proved there  
can be extended to  
the Mori case by using different technical tools. 
In Theorem \ref{wMori} we characterize $w$-divisorial  
Mori domains and in Theorems \ref{wvMori} and  \ref{Mori2}  
we study  $w$-divisoriality of their overrings. 
In particular, we show that generalized rings of fractions of $w$-divisorial  
Mori domains are $w$-divisorial and we prove that a domain 
whose $t$-linked overings are all $w$-divisorial is Mori if and only if  
it has $t$-dimension one. 

We thank the referee for his/her careful reading and relevant observations.

\smallskip 
Throughout this paper $R$ will denote an integral domain 
with quotient field $K$ and we will assume that $R \neq K$.  
 
We shall use the language of star-operations. A {\it star operation} is  
a map $I\to I^*$  from the set $F(R)$ of nonzero fractional ideals of $R$ to 
itself such that:  
 
(1) $R^* = R$ and $(aI)^* = aI^*$, for all $a \in K \smallsetminus 
\{0\}$; 
 
(2) $I \sub I^*$ and $I \sub J \ra I^* \sub J^*$; 
 
(3) $I^{**} = I^*$. 
 
General references for systems of ideals and star operations are  
 \cite{g1, Gr, HK, J}. 
 
 A star operation $*$ is of {\it finite type} if 
 $I^* = \cup \{J^*\, ;\, J \sub I$  and $J$ is fini\-te\-ly 
ge\-ne\-ra\-ted\}, 
 for each $I \in F(R)$. 
To any  star operation $*$, we can associate a star operation 
$*_{f}$ of finite type 
by defining $I^{*_{f}}= \cup J^*$, with the 
union taken over all finitely generated ideals $J$ contained 
in $I$. Clearly $I^{*_{f}} \sub I^*$.  
A nonzero ideal $I$ is $*$-{\it finite} if $ I^* =J^*$ for  
 some finitely generated ideal $J$. 
 
The identity is a star operation, called the $d$-operation. 
The $v$- and the 
$t$-operations are the best known nontrivial star operations and are 
defined in the following way. 
For a pair of nonzero ideals $I$ 
and $J$ of a domain $R$ we let $(J\colon I)$ denote the set $\{x\in K  
\, ; \, xI\subseteq J\}$.  
We set $I_v=(R\colon (R\colon I))$ and $I_t=\cup J_v$ with the 
union taken over all finitely generated ideals $J$ contained 
in $I$. Thus the $t$-operation is the finite type  
star operation associated to the $v$-operation. 
 
A nonzero fractional ideal $I$ is called a $*$-{\it ideal} if $I= I^*$.  
If $I = I_{v}€$  we say that $I$ is {\it divisorial}.  
For each star operation $*$, we have $I^* \sub  
I_{v}€$, thus each divisorial ideal is a $*$-ideal. 
 
The  set $F_{*}(R)$ of $*$-ideals of  
 $R$ is a semigroup with respect to the  
 $*$-{\it multiplication}, defined by $(I,J) \to (IJ)^*$, with unity $R$.  
 We say that an ideal $I \in F(R)$ is $*$-{\it invertible} if $I^*$ is a  
 unit in the  
 semigroup $F_{*}(R)$. In this case the $*$-{\it inverse} of  
$I$ is $(R:I)$. 
 Thus $I$ is  $*$-invertible if and only if $(I(R:I))^*=R$. Invertible  
 ideals are ($*$-invertible) $*$-ideals. 
 
A prime $*$-ideal is also called a $*$-{\it prime}.  
A $*$-{\it maximal} ideal is an ideal 
that is maximal in the set of the proper  
$*$-ideals. A $*$-maximal ideal  (if it exists) is a prime ideal. 
If $*$ is a star operation of finite type, an easy application of  
Zorn's Lemma shows that the set $*$-Max$(R)$ of the $*$-maximal ideals of 
$R$ is not empty. Moreover, for each $I\in F(R)$, $I^* = \cap_{M\in 
*\op{-Max}(R)}I^*R_{M}$; in particular  
$R = \cap_{M\in *\op{-Max}(R)}R_{M}$ \cite{Gr}. 
 
The $w$-operation is the star operation  defined by  
setting $I_{w} = \cap_{M\in t\text{-Max}(R)}IR_{M}$.  
 An equivalent definition is obtained by setting 
 $I_{w}€ = \cup \{(I\colon  
J) \; ; \; J$ is finitely generated and $(R\colon J) = R\}$. 
By using the latter definition, one can see that the notion of $w$-ideal coincides  
with the notion of {\it semi-divisorial} ideal introduced by S. Glaz  
and W. Vasconcelos in 1977 \cite{GV}. As a star-operation, the $w$-operation  
was first considered by E. Hedstrom and E. Houston in  
$1980$ under the name of $F_{\infty}€$-operation \cite{HH}.  
Since 1997 this star operation was intensely studied by  Wang Fanggui and  
R. McCasland in a more general context. In particular they showed that  
the notion of $w$-closure is a very useful tool in the study of Strong 
Mori domains \cite{FM, FM2}. 
 
The $w$-operation is of finite type. 
We have $w$-Max$(R) = \tmax(R)$ and $IR_{M}€= I_{w}€R_{M}€\sub  
I_{t}€R_{M}€$, for 
each $I \in F(R)$ and $M \in \tmax(R)$. Thus $I_w \sub I_t \sub I_v$. 
 
We denote by $\tspec(R)$  
the set of $t$-prime ideals of $R$. Each height one prime is a $t$-prime  
and each prime minimal over a $t$-ideal is a $t$-prime. 
We say that $R$ has $t$-{\it dimension  
one} if each $t$-prime ideal has height one.

%%%%%%%%%%%%% w-DIVISORIAL DOMAINS %%%%%%%%%%%%%%%%%%%%%%%%%%%%%%%%%% 

\section{$w$-divisorial domains} 
 
A {\it divisorial domain} is a domain such  that each ideal is  
divisorial  \cite{BS} and we say that a domain $R$ is {\it $w$-divisorial}  
if each $w$-ideal is divisorial, that is $w = v$.  Since $I_{w}\sub I_{t}  
\sub I_{v}$, for each nonzero fractional ideal $I$, then $R$ is 
$w$-divisorial if and only if $w = t =v$.  A  domain with the 
property that $t=v$ is called in \cite{HZ} a {\it  
$TV$-domain}. Mori domains (i. e. domains satisfying the ascending chain  
condition on proper divisorial ideals) are $TV$-domains.  
A domain such that $w=t$ is called a {\it $TW$-domain} \cite{M2}.  
An important class of $TW$-domains is the class  
of $PvMD$s; in fact a $PvMD$ is precisely an integrally closed  
$TW$-domain \cite[Theorem 3.1]{K}. (Recall that a domain  
$R$ is a {\it  Pr\"ufer $v$-multiplication domain}, for short a $PvMD$,  
if $R_{M}€$ is a valuation domain for each $t$-maximal ideal $M$ of  
$R$.)  Since a Krull domain is a Mori  
$PvMD$, a Krull domain is a $w$-divisorial domain.  
An example due to M. Zafrullah  
shows that in general  $w\neq t \neq v$ \cite[Proposition 1.2]{M2}.  
Also there exist $TV$-domains and $TW$-domains that are not  
$w$-divisorial \cite[Example  2.7]{M2}.  
 
If $R$ is a  Pr\"ufer  domain, in particular a valuation domain, then 
$w$-divisoriality coincides with divisoriality, because each ideal 
of a  Pr\"ufer  domain is a $t$-ideal.

\begin{proposition}\label{P1} 
A $w$-divisorial domain $R$ is divisorial if and only if each maximal  
ideal of $R$ is a $t$-ideal. Hence a one-dimensional $w$-divisorial  
domain is divisorial.  
\end{proposition}  
 
\begin{proof}If  
each maximal ideal of $R$ is a $t$-ideal, then each ideal of $R$ is a  
$w$-ideal by \cite[Proposition 1.3]{M2}. Hence, if  $R$ is  
$w$-divisorial it is also divisorial. The converse is clear.  
\end{proof}

Following \cite{AZ}, we say that  
a nonempty family $\La$ of nonzero prime ideals of $R$ is of {\it finite  
character}  
if each nonzero element of $R$ belongs to  
at most finitely many members of $\La$ and we say that $\La$ 
is {\it independent} if no two members of $\La$  contain a common  
nonzero prime ideal. We observe that a family of primes is independent  
if and only if no two members  
of $\La$  contain a common $t$-prime ideal. In fact a minimal prime of  
a nonzero principal ideal is a $t$-ideal.  
 
The domain $R$ has finite character  
(resp., $t$-finite character) if  $\Max(R)$ (resp., $\tmax(R)$) is of  
finite character. 
If the set $\Max(R)$ is independent of finite  
character, the domain $R$ is called by E. Matlis  an {\it $h$-local}  
domain  
\cite{Ma}; thus $R$ is $h$-local if it has finite character and  
each nonzero prime ideal is contained  
in a unique maximal ideal. 
A domain $R$ such that $\tmax(R)$ is independent of finite  
character is called in \cite{AZ} a {\it weakly Matlis domain}; hence 
 $R$ is a  weakly Matlis  
domain  if it has $t$-finite character and each $t$-prime ideal is contained  
in a unique $t$-maximal ideal. 
 
 Clearly, a domain of $t$-dimension one is a weakly Matlis domain if and 
only if it  has $t$-finite character. A one-dimensional domain is a 
weakly Matlis domain if and only if it is $h$-local;  if and only if it 
has finite character. 
 
We recall that any  
$TV$-domain, hence any $w$-divisorial domain, has $t$-finite character  
by \cite[Theorem 1.3]{HZ}. The  main result of this section shows  
that $w$-divisorial domains form a distinguished class of weakly Matlis  
domains.   
 
We start by proving some technical properties of weakly  
Matlis domains.  
 
\begin{lemma}\label{L1}  Let $R$ be an integral domain. The  
following conditions are equivalent:  \begin{itemize}  
\item[(1)] $R$ is a weakly  
Matlis domain;   
\item[(2)] For each $t$-maximal ideal $M$ of $R$ and  
a collection $\{I_\al\}$ of $w$-ideals of $R$ such that $\cap_\al  
I_\al\neq 0$, if $\cap_\al I_\al\sub M$, then  $I_\al\sub M$ for some  
$\al$.  \end{itemize}\end{lemma}  
 
\begin{proof} $(1)\ra(2)$ follows from \cite[Corollary 4.4 and  
Proposition 4.7]{AZ}, by taking ${\mc F}=\tmax(R)$ and then  
$\ast_{\mc F}=w$.  
 
$(2)\ra(1)$. First, we show that each $t$-prime ideal is contained in  
a unique $t$-maximal ideal. We adapt the proof of \cite[Theorem  
2.4]{H}.  
Let $P$ be a $t$-prime which is contained in two distinct  
$t$-maximal ideals $M_1$ and $M_2$. Let $\{I_\al\}$ be the set of all  
$w$-ideals of $R$ which contain $P$ but are not contained in  
$M_1$. Such a collection is nonempty since $M_{2}€$ is in it. Let  
$I=\cap I_\al$. Then $I\nsub  M_1$ and $I \sub M_{2}€$.  
Take $x\in I\sm M_1$. Since $x^2\notin M_1$, then $(P+x^2R)_w\in  
\{I_\al\}$ and so $x\in  (P+x^2R)_w$. Thus $x\in  
(P+x^2R)R_{M_2}\neq R_{M_2}$ and $sx=p+x^2r$ for some $s\in R\sm  
M_2$, 
$p \in P$ and $r\in R$.  
Whence $(s-rx)x=p\in P \sub M_{1}€\cap M_{2}€$. Now $s-rx\notin P$ 
because $s\notin M_2$  and $rx \in I \sub M_{2}$. But also $x\notin  
P$, since $x\notin  M_1$; a contradiction because $P$ is prime.   
 
Next we show that $R$ has  
$t$-finite character. Let $0\neq x \in R$ and $\{M_\be\}$ be the set of  
all $t$-maximal ideals of $R$ which contain $x$. For a fixed $\be$, let  
$A_\be$ be the intersection of all $w$-ideals of $R$ which  
contain $x$ but are not contained in $M_\be$. By assumption  
$A_\be\nsub  M_\be$. Set $A=\sum_\be A_\be$. Then $x\in A$ and $A$ is  
contained in no $M_\be$. Hence $A_t=R$. Let $F=(a_{\be_1},  
a_{\be_2},\ldots  , a_{\be_n})$, where $a_{\be_i}\in A_{\be_i}$, be a  
finitely generated ideal of $R$ such that $F_t=R$. Now, if $M_\be\notin  
\{M_{\be_1}, M_{\be_2},\ldots  , M_{\be_n}\}$, necessarily  
$M_\be\supseteq F$, which is impossible because $M_\be$ is a proper  
$t$-ideal and $F_t=R$. We conclude that 
$\{M_\be\}=\{M_{\be_1}, M_{\be_2},\ldots,M_{\be_n}\}$ is finite.  
\end{proof} 
 
 \begin{lemma}\label{L2}  Let $R$ be a $w$-divisorial domain, 
$M$  a $t$-maximal ideal of $R$ and $\{I_\al\}$ a collection of  
$w$-ideals of $R$ such that $\cap_\al I_\al\neq 0$. If $\cap_\al  
I_\al\sub M$, then  $I_\al\sub M$ for some $\al$. \end{lemma}  
 
\begin{proof} Set $A= \cap_\al I_\al$. Since $R$  
is a $TW$-domain, then the $I_\al$'s and $A$ are $t$-ideals. Since $R$  
is also a $TV$-domain, by \cite[Lemma 1.2]{HZ}, if $I_{\al} \nsubseteq  
M$, for each $\al$, then $A \nsubseteq M$. 
\end{proof}  
 
\begin{lemma}\label{L3}  If $R$ is a weakly Matlis domain,  then  
$I_vR_M =(IR_M)_v$, for each nonzero fractional ideal $I$ and each  
$t$-maximal ideal  $M$. \end{lemma}  
 
\begin{proof} Apply \cite[Corollary 5.3]{AZ} for ${\mc  
F}=\tmax(R)$. \end{proof} 
 
  We are now ready to prove the $t$-analogue of \cite[Proposition  
5.4]{BS}, which states that a domain $R$ is divisorial if and only if  
it is $h$-local and $R_M$ is a divisorial domain, for each maximal  
ideal $M$. Local divisorial domains have been studied in \cite[Section  
5]{BS} and completely characterized in \cite[Section 2]{B2}.  
 
\begin{theorem}\label{Loc}Let $R$ be an integral domain.  
The following conditions are equivalent:  
\begin{itemize}  
     \item[(1)] $R$ is a  
$w$-divisorial domain;   
\item[(2)] $R$ is a weakly Matlis domain and  
$R_M$ is a divisorial domain, for each $t$-maximal ideal  $M$;  
\item[(3)] $R$ is a $TV$-domain and $R_M$ is a divisorial domain, for  
each 
$t$-maximal ideal  $M$;  
\item[(4)] $IR_{M}=(IR_M)_v=I_vR_M$, for each  
nonzero fractional ideal $I$ and each $t$-maximal ideal  $M$.  
\end{itemize}\end{theorem}  
 
\begin{proof} $(1)\ra(2)$. That $R$ is a weakly Matlis domain  
follows from Lemmas \ref{L2} and \ref{L1}. Now let $M$ be a $t$-maximal  
ideal of $R$ and $I = JR_M$ a nonzero ideal of $R_M$, where  
$J$ is an ideal of $R$. By Lemma \ref{L3}, we have  
$I_v=(JR_M)_v=J_vR_M$.  
Since $J_v=J_w$, then $I_v=J_wR_M=JR_M=I$. Hence $R_M$ is a  
divisorial domain. 
 
$(2)\ra(4)$ follows from Lemma \ref{L3}.  
 
 $(4)\ra(1)$. Let $I$ be a nonzero fractional ideal of $R$.  Then  
$I_w=\cap_{M\in \tmax (R)} IR_M = 
\cap_{M\in \tmax (R)} I_vR_M=I_v$.  
Whence $R$ is a  $w$-divisorial domain.  
 
$(1)\ra(3)$ via (2). 
 
$(3)\ra(4)$.  Since $t=v$ in $R$ and $d=t=v$ in $R_{M}Û$, 
 for each nonzero fractional ideal $I$ and each $t$-maximal ideal 
$M$ of $R$, we have 
 
$IR_{M}Û=(IR_{M})_{v}Û=(IR_{M})_{t}Û=(I_{t}ÛR_{M})_{t}=I_{t}ÛR_{M}=I_{v}ÛR_{M}$. 
\end{proof} 
 
Any almost Dedekind domain that is not Dedekind provides an example  
of a locally divisorial domain that is not $w$-divisorial, because it  
is not of finite character \cite[Theorem 37.2]{g1}.  
 
\begin{corollary}\label{tdim1} Let $R$ be a domain of $t$-dimension  
one. Then $R$ is $w$-divisorial if and only if $R$ has $t$-finite  
character and  
$R_{P}€$ is divisorial, for each height one prime $P$. 
\end{corollary} 
 
%%%%%%%%%%%%%%%%%%%%%%%%%%%%%%%%%%%%% LOCALIZATIONS  
%%%%%%%%%%%%%%%%%%%%%%%%%%%%%%%%%%%%%%%% 

 \section{Localizations of $w$-divisorial domains}  
 
A domain whose overrings are all divisorial is called {\it totally  
divisorial} \cite{BS}. Not all divisorial domains are totally  
divisorial \cite[Remark 5.4]{H}; in fact a valuation domain $R$ is  
divisorial if and only if its maximal ideal is principal \cite[Lemma  
5.2]{H}, but it is totally divisorial if and only if it is strongly  
discrete \cite[Proposition 7.6]{BS}, equivalently $PR_{P}€$ is a  
principal  
ideal for each prime ideal $P$ of $R$ \cite[Proposition 5.3.8]{fhp}.  
Since for valuation domains divisoriality coincides with  
$w$-divisoriality and each overring of a valuation domain is a  
localization at a certain ($t$-)prime, 
 we see that $w$-divisoriality is not stable under localization at  
$t$-primes. 
 
 We say that an integral domain $R$ is a {\it strongly  
$w$-divisorial domain} (resp.,  a {\it strongly  divisorial domain}) 
if $R$ is $w$-divisorial (resp., divisorial) and $R_P$ is a divisorial  
domain for each $P\in\tspec(R)$ (resp., $P\in\spec(R)$).  
Note that if $R$ is strongly $w$-divisorial (resp., strongly  
divisorial),  
then $R_P$ is strongly divisorial for each $P\in\tspec(R)$ (resp., for  
each $P\in\spec(R)$). 
 
By Theorem  
\ref{Loc} (resp.,  \cite[Proposition 5.4]{BS}), $R$ is a strongly  
$w$-divisorial domain (resp.,  a strongly  divisorial domain) 
if and only if $R$ is a weakly Matlis domain (resp., an $h$-local  
domain)  
and $R_P$ is a divisorial  
domain for each $P\in\tspec(R)$ (resp., $P\in\spec(R)$).  
 
If $R$ has $t$-dimension one,  then $R$ is $w$-divisorial if and only  
if it is strongly  
$w$-divisorial. 
\medskip 
 
In this section we shall study the extension of $w$-divisoriality and  
divisoriality 
 to distinguished classes of generalized rings of fractions such as  
localizations at ($t$-)prime 
ideals, ($t$-)flat overrings and ($t$-)subintersections. 
 
We recall the requisite definitions. A  
nonempty family  
$\mc F$ of nonzero ideals of a domain $R$ is said to be a {\it  
multiplicative system} 
of ideals  if $IJ \in \mc F$, for each 
$I, J \in \mc F$.  If $\mc F$ is a multiplicative system, the 
set of ideals of $R$ containing some ideal of $\mc F$ is  still 
a  multiplicative system, which is  called the {\it saturation  
of $\mc F$} and is  denoted by Sat$(\mc F)$. A multiplicative system 
$\mc F$ is said to be {\it saturated} if $\mc F =$ Sat$(\mc F)$. 
 
If $\mc F$ is a multiplicative system of ideals,  the overring $R_{\mc  
F}:=\cup\{(R:J);\, J\in {\mc F}\}$  of $R$ is called the {\it  
generalized  
ring of fractions} of $R$ with respect to $\mc F$. For any fractional  
ideal  
$I$ of $R$, $I_{\mc F}:=\cup\{(I:J);\, J\in {\mc  
F}\}$ is a fractional ideal of $R_{\mc F}$ and $IR_{\mc F}\sub I_{\mc  
F}$. Clearly $I_{\mc F}€ = I_{\op{Sat}(\mc F)}€$. 
 
The map $P\mapsto P_{\mc F}$ is an 
order-preserving bijection between the set of prime ideals $P$ of $R$ such that 
$P\notin \op{Sat}({\mc F})$ and the set of prime ideals $Q$ of  
$R_{\mc F}$ such 
that $JR_{\mc F}\not\subseteq Q$ for any $J\in {\mc F}$, with inverse map 
$Q\mapsto Q\cap R$. In addition, $R_P=(R_{\mc F})_{P_{\mc F}}$ for each prime 
ideal $P\notin \op{Sat}({\mc F})$.  
If $Q$ is a $t$-prime ideal of $R_{\mc F}$, 
then $Q\cap R$ is a $t$-prime ideal of $R$ \cite[Proposition 1.3]{gab2}. 
 
If $\La$ is a nonempty  
family of nonzero prime ideals of $R$, the set ${\mc F}(\La)=\{ J \,;  
\; J \sub R 
\; $is an ideal and $\; J\nsub P \; $for  
each $P\in \La\}$ is a saturated  
multiplicative system of ideals and  $I_{{\mc  
F}(\La)}=\cap\{IR_P;\, P\in\La\}$, for each fractional ideal $I$ of  
$R$; in particular $R_{{\mc F}(\La)}=\cap\{R_P;\,  
P\in\La\}$.  
A generalized ring of fractions of type $R_{{\mc F}(\La)}$ is called a  
{\it subintersection of $R$}; when 
$\La\sub\tspec(R)$, we say that $R_{{\mc F}(\La)}$ is a 
{\it $t$-subintersection} of $R$.  
 
A multiplicative system of ideals ${\mc F}$ of $R$ is {\it finitely  
generated} if each ideal $I\in {\mc F}$ contains a finitely generated  
ideal $J$ which is still in ${\mc F}$. As in \cite{gab2}, 
 we say that ${\mc F}$ is a {\it $v$-finite} multiplicative system  
if each $t$-ideal  $I\in \op{Sat}({\mc F})$ contains a  
finitely generated ideal $J$ such that $J_v\in\op{Sat}({\mc F})$.  
A finitely generated multiplicative system is $v$-finite. 
If ${\mc F}$ is $v$-finite, 
the set $\La$  of $t$-ideals which are maximal with respect to  
the property of not being in  
Sat$(\mc F)$ is not empty, $\La \sub t\op{-Spec}(R)$, ${\mc  
F}(\La)$ is $v$-finite and $T =R_{{\mc F}(\La)}$ \cite[Proposition 1.9  
(a) and (b)]{gab2}. 
 
An overring  $T$  of $R$ is said to be {\it $t$-flat} over  
$R$  if $T_M=R_{M \cap R}$, for each $t$-maximal ideal $M$ of $T$  
\cite{KP}, 
equivalently $T_Q=R_{Q \cap R}$, for each $t$-prime ideal $Q$ of $T$  
\cite[Proposition 2.6]{elb1}.  
Flatness implies $t$-flatness, but the converse is not true  
\cite[Remark 2.12]{KP}.  
By \cite[Theorem 2.6]{elb1}, $T$ is $t$-flat over $R$ if and only if there  
exists a $v$-finite multiplicative system $\mc F$ of $R$ such that $T =  
R_{\mc F}€$. Thus $T$ is $t$-flat if and only if $T =R_{{\mc F}(\La)}$,  
where $\La$ is a family of pairwise incomparable  
$t$-primes of $R$ and ${\mc F}(\La)$ is $v$-finite. 
It follows that a $t$-flat overring of $R$ is a 
 $t$-subintersection of $R$.  
 
In turn, any 
 generalized ring of fractions is a  
$t$-linked overring; but the converse does not hold in general  
\cite[Proposition 2.2]{DHLZ}. 
We recall that an overring $T$ of an integral domain $R$ is  
$t$-{\it linked} over $R$ 
if, for each nonzero finitely generated ideal $J$ of $R$ such that  
$(R:J)=R$, we have 
$(T:JT)=T$ \cite {DHLZ}. This is equivalent to say that $T  
=\cap T_{R\setminus P}$, where $P$ ranges over the $t$-primes of $R$  
\cite[Proposition 2.13(a)]{DHLZ}.  
\medskip 
 
It is well known that if $P$ is a $t$-prime ideal of $R$, then $PR_{P}€$ need not  
 be a $t$-ideal of $R_{P}€$. When $PR_{P}€$ is a  
$t$-prime ideal, $P$ is called by M. Zafrullah a {\it well behaved} $t$-prime 
\cite[page 436]{Z}.  
We prefer to say that $P$ $t$-{\it localizes} or that it is a 
$t$-{\it localizing prime}.  Height-one prime ideals  
and divisorial $t$-maximal primes, e. g. $t$-invertible $t$-primes, are  
examples of $t$-localizing primes. 
 
A large class of domains with the property that  
each  $t$-prime ideal $t$-localizes is the class of $v$-coherent domains. 
We recall that a domain $R$ is called $v$-{\it coherent} if the ideal  
$(R:J)$ is $v$-finite whenever $J$ is finitely generated. This class of  
domains  properly includes $PvMD$'s, Mori domains  
and coherent domains \cite{N, GH1}. 
 
If $R$ is a $w$-divisorial (resp., strongly $w$-divisorial) domain,  
then each $t$-maximal (resp., $t$-prime) ideal $t$-localizes. 
 
\begin{lemma}\label{max} Let $\La$ be a set of  
$t$-localizing $t$-primes of $R$. Then: 
     \begin{itemize} 
\item[(1)] $P_{{\mc F}(\La)}€\in t\op{-Spec}(R_{{\mc F}(\La)})$, for  
each $P\in \La$. 
\item[(2)] If  ${\mc F}(\La)$ is  $v$-finite,  $t\op{-Max}(R_{{\mc  
F}(\La)}) = \{P_{{\mc F}(\La)}€\; ; P$ maximal in $\La\}$. 
\end{itemize} 
\end{lemma} 
\begin{proof} Set $\mc F ={\mc F}(\La)$ and $T=R_{\mc F}€$.  
     
(1). Let $P\in\La$. Since $R_{P} = T_{P_{\mc F}}$ and by hypothesis 
         $PR_P=P_{\mc F}€T_{P_{\mc F}}$  
is a $t$-ideal, then $P_{\mc F}=P_{\mc F}€T_{P_{\mc F}}\cap  
T$ is a $t$-ideal of $T$.  
 
(2).     Since $P_{\mc F}$ is a $t$-ideal by part (1), 
         we can apply \cite[Proposition 1.9 (c)]{gab2}.  
\end{proof}

\begin{proposition}\label{WMover} Let $\La$ be a set of pairwise incomparable  
$t$-localizing $t$-primes of $R$. Then: 
     \begin{itemize} 
         \item[(1)] $\La$ is independent of finite character if and only if  
${\mc F}(\La)$ is $v$-finite and $R_{{\mc F}(\La)}$ is a weakly Matlis  
domain. 
\item[(2)] If $R_{{\mc F}(\La)}$ is $w$-divisorial, then $\La$ is  
independent  
of finite character.  
\end{itemize} 
     \end{proposition} 
     
     \begin{proof}     Set $\mc F = {\mc F}(\La)$ and $T = R_{\mc F}€$.  
     
         (1).     If ${\mc F}$ is $v$-finite, by Lemma  
         \ref{max}(2) we have 
$\tmax(T) = \{P_{\mc F} \; ; P\in \La\}$. It follows  
that $\La$ is independent of finite character if and only if  
$\tmax(T) = \{P_{\mc F} \, ; P\in \La\}$  
is independent of finite character, that is $T$ is a weakly Matlis  
domain.  
On the other hand, if $\La$ is of finite character, then ${\mc F}$ is  
$v$-finite by  
\cite[Lemma 1.16]{gab2}. 
     
     (2). Since $T$ is a weakly Matlis  
     domain, by part (1) it suffices to show that $\La$ is of finite  
     character. 
         
     By Lemma \ref{max}(1), $P_{\mc F}€$ is a $t$-prime of $T$, for each $P  
\in \La$.  
We show that each proper divisorial ideal of $T$ is contained in some 
     $P_{\mc F}€$. 
     We have $T = \cap_{P\in \La}€ R_{P}€= \cap_{P\in \La}€ T_{P_{\mc  
F}€}€$.  
     If $I$ is a proper divisorial ideal of $T$, there is $x  
     \in K \sm T$ (where $K$ is the quotient field of $R$)  
         such that $I \sub x^{-1}T\cap T$. Since $x \notin T$,  
     there exists $P \in \La$ such that  $x \notin T_{P_{\mc F}€}€$,  
     equivalently $x^{-1}T\cap T \sub P_{\mc F}$. 
     
     Since $t=v$ on $T$, we conclude that $\tmax(T) = 
     \{P_{\mc F}€\; ; P\in \La\}.$ Since  $T$ has $t$-finite  
     character, it follows that $\La$ is of finite character.  
\end{proof} 
 
\begin{theorem}\label{tmaxtsub} Let $R$ be a $w$-divisorial domain. If  
$\La \sub  
t\op{-Max}(R)$, then $R_{{\mc F}(\La)}€$ is a $t$-flat $w$-divisorial  
overring of $R$. 
\end{theorem} 
 
\begin{proof} Since $R$ is a weakly Matlis domain (Theorem \ref{Loc}),  
     $\tmax(R)$ is independent of finite character; thus $\La$ has the  
same properties. In  
addition, each $t$-maximal ideal is a $t$-localizing prime ideal.  
It follows that ${\mc F}(\La)$ is $v$-finite and $T := R_{{\mc  
F}(\La)}€$ is a $t$-flat weakly Matlis domain (Proposition  
\ref{WMover}(1)).  
By Lemma \ref{max}(2), for each $N \in t\op{-Max}(T)$, 
there exists $M\in \La$ such that $N  = M_{{\mc F}(\La)}$. 
It follows that $T_{N}= R_{M}$ is divisorial and so  
 $T$ is $w$-divisorial by Theorem \ref{Loc}. 
\end{proof} 
 
As we have mentioned above, the localization of a $w$-divisorial domain  
at a $t$-prime need not be a ($w$-)divisorial domain. Thus  
Theorem \ref{tmaxtsub} does not  
hold for an arbitrary $\La \sub t\op{-Spec}(R)$. However,  
under the hypothesis that $R$ is strongly 
 $w$-divisorial, we have a satisfying result.

\begin{theorem} \label{subint} Let $R$ be a strongly $w$-divisorial  
domain and  
     $\La$ a set of pairwise incomparable $t$-primes of $R$. The following  
conditions  
     are equivalent: 
     \begin{itemize} 
         \item[(1)] $R_{{\mc F}(\La)}$ is $w$-divisorial;  
         \item[(2)] $R_{{\mc F}(\La)}$ is strongly $w$-divisorial;  
        \item[(3)] $R_{{\mc F}(\La)}$ is a $t$-flat weakly Matlis domain; 
\item[(4)] $R_{{\mc F}(\La)}$ is a $t$-flat TV-domain; 
 \item[(5)] $\La$ is independent of finite character.  
     \end{itemize} 
     \end{theorem} 
     \begin{proof}  
     Set $\mc F ={\mc F}(\La)$ and  $T = R_{\mc F}$. Since $R$ is  
strongly $w$-divisorial, each $P\in \La$ $t$-localizes. 
     
     $(1)\ra(5)$ by Proposition \ref{WMover}(2). 
 
$(5)\ra(3)$. By Proposition \ref{WMover}(1).  
 
$(3)\ra(2)$. If $Q$ is a $t$-prime of $T$, then $P = Q\cap R \in  
\tspec(R)$ and 
  $T_{Q}€=R_{P}€$ is divisorial. Whence $T$ is strongly  
  $w$-divisorial.

$(3)\lra(4)$ By $t$-flatness, $T_M$ is divisorial for each $t$-maximal ideal 
 $M$. Thus we can apply Theorem \ref{Loc}. 
 
$(2)\ra(1)$ is obvious. 
\end{proof} 
 
Divisorial flat overrings of a strongly divisorial domain have a 
 similar characterization.   Recall that an overring $T$ of $R$ is flat if $T_M=R_{M\cap R}$, for each maximal ideal $M$ of $T$; in this case $T= R_{{\mc F}(\La)}$, where $\La$ is a set of pairwise incomparable prime ideals of $R$. 
 
\begin{corollary}\label{flat} Let $R$ be a strongly divisorial  domain and $T= R_{{\mc F}(\La)}$ a flat overring, where $\La$ is a set of pairwise incomparable prime ideals of $R$. The following  
conditions are equivalent: 
\begin{itemize} 
 
\item[(1)] $T$ is divisorial; 
 
\item[(2)] $T$ is strongly divisorial; 

\item[(3)] $T$ is $h$-local;
 
\item[(4)]  $\La$ is independent of finite character.\end{itemize} 
\end{corollary} 
\begin{proof}  $(1)\lra (3)$. By \cite[Proposition 5.4]{BS}, $T$ is divisorial if and only if it is $h$-local and locally divisorial. But, since $T$ is flat and $R$ is strongly divisorial, for each maximal ideal $M$ of $T$, $T_M= R_{M\cap R}$ is divisorial.

 $(1) \ra (2)$.  Since $T$ is flat and $R$ is strongly divisorial, then $T_Q=R_{Q\cap R}$ is divisorial, for each prime ideal $Q$ of $T$.

     $(2) \ra (4)$. 
     Since $R$ and $T$ are divisorial, then $d = w= t= v$ in $R$ and $T$. Thus we can apply  
    Theorem \ref{subint} ($(2) \ra (5)$). 
     
     $(4)\ra (1)$.  Since $d = w= t= v$ in $R$, by  Theorem \ref{subint} ($(5) \ra (1)$),  $T$ is $w$-divisorial. To prove that $T$ is  
     divisorial, we show that each maximal  
     ideal of $T$ is a $t$-ideal (Proposition \ref{P1}). If $M$ is a maximal ideal of $T$, by  
     flatness we have $T_{M}= R_{M\cap R}$. Since $R$ is strongly divisorial, 
     $MT_{M}$ is a $t$-ideal and so $M = MT_{M}\cap T$ is a $t$-ideal. 
     \end{proof} 
  
\begin{corollary} \label{stflat} 
     Let $R$ be an integral domain.  
The following conditions are equivalent:\begin{itemize} 
\item[(1)] Each $t$-flat overring of $R$ is strongly $w$-divisorial; 
\item[(2)] $R$ is strongly $w$-divisorial and each $t$-flat overring is a weakly Matlis domain; 
\item[(3)] $R$ is strongly $w$-divisorial and each $t$-flat overring is a $TV$-domain;  
\item[(4)] $R$ is strongly $w$-divisorial and each family  
$\La$ of pairwise incomparable $t$-primes of $R$ such that ${\mc  
F}(\La)$ is $v$-finite is independent of finite  
character.\end{itemize}  
\end{corollary}  
\begin{proof} By Theorem \ref{subint}, recalling that an overring  
$T$ is $t$-flat over $R$ if and only if $T =R_{{\mc F}(\La)}$,  
where $\La$ is a family of pairwise incomparable  
$t$-primes of $R$ and ${\mc F}(\La)$ is $v$-finite. 
\end{proof} 
 
In order to study $t$-subintersections, we need the following technical lemma. 
 
\begin{lemma} \label{acc} 
Let $R$ be an integral domain and $\mc C$ an ascending  
chain of $t$-localizing $t$-primes of $R$. If $R_{\mc F(\mc C)}$ is a $TV$-domain, then  
${\mc C}$ is stationary. 
\end{lemma}  
\begin{proof} 
 Let ${\mc C}=\{P_{\alpha}\}$ and set ${\mc F}={\mc F(\mc C)}$ and $T=R_{\mc F}$. By Lemma \ref{max}(1),  
$(P_{\alpha})_{\mc F}$ is a $t$-prime ideal of $T$, for each $\alpha$.  
It follows that 
$M=\cup_{\alpha}(P_{\alpha})_{\mc F}$ is a proper $t$-prime ideal of $T$  
(since it is an  
ascending union of $t$-primes) and so $M$ is divisorial (because $T$ is a  
$TV$-domain). 
We have  
$T=\cap_{\alpha}T_{R\sm {P_{\alpha}}}$; thus the map $I\mapsto  
I^\star=\cap_{\alpha}IT_{R\sm {P_{\alpha}}}$ defines a star operation 
 on $T$. Since $M$ is divisorial, we have 
$M^\star\sub M$; so that $M^\star$ is a proper ideal. It  
follows 
that there exists $\alpha$ such that $M\cap R\sub P_\alpha$. Hence  
$M\cap R= P_\alpha$ and so $P_\beta = P_\alpha$ for $\beta \geq \alpha$. 
\end{proof} 
 
 \begin{theorem}\label{stsubint}  
     Let $R$ be an integral domain.  
The following conditions are equivalent:\begin{itemize} 
 
\item[(1)] Each $t$-subintersection of $R$ is strongly $w$-divisorial; 
 
\item[(2)] $R$ is a strongly $w$-divisorial domain which satisfies the  
ascending chain condition on  
$t$-prime ideals and each family  
$\La$ of pairwise incomparable $t$-primes of $R$ is independent  
of finite character. 
\end{itemize}  
\end{theorem} 
\begin{proof} $(1)\ra(2)$. Clearly $R$ is a strongly $w$-divisorial  
domain. If  
$\La$ is a set  of pairwise incomparable $t$-prime ideals, then by  
assumption  
$R_{\mc F(\La)}$ is strongly $w$-divisorial. Hence $\La$ is independent  
of finite character, by Theorem \ref{subint}. It remains to show that  
$R$ has the ascending chain condition on $t$-prime ideals. This follows from  Lemma \ref{acc}.  
In fact, if 
$\mc C $ is an ascending  
chain of $t$-prime ideals of $R$, $R_{\mc F(\mc C)}$ is strongly $w$-divisorial. 
Hence each $t$-prime in $\mc C$ $t$-localizes and it follows that $\mc C$ is stationary. 
 
$(2)\ra(1)$. Let $R_{{\mc F}(\La)}$  be a $t$-subintersection of $R$. 
 By the ascending chain condition on $t$-prime ideals,  
 $\La$ has maximal elements; thus we  
can assume that $\La$ is a set of  
pairwise incomparable $t$-primes. The conclusion follows from Theorem \ref{subint}.  
\end{proof}

\begin{corollary} Let $R$ be a domain. 
If each $t$-subintersection of $R$ is strongly $w$-divisorial, then each  
$t$-subintersection of $R$ is $t$-flat. 
\end{corollary}  
\begin{proof} If each $t$-subintersection of $R$ is strongly $w$-divisorial,  
     then $R$ satisfies 
the ascending chain condition on $t$-primes (Theorem \ref{stsubint}).  
Thus each $t$-subintersection is of type 
$R_{\mc F(\La)}$, where $\La$ is a family of pairwise incomparable $t$-primes. 
By Theorem \ref{subint}, $R_{\mc F(\La)}$ is $t$-flat. 
\end{proof} 
 
\begin{remark}\label{subint2}  \rm If each subintersection of the domain $R$ 
     is strongly divisorial, then clearly  
$R$ is  strongly divisorial. In addition, since $d=w=t=v$ on $R$, then $R$ 
satisfies the  ascending chain condition on  
prime ideals and  each family $\La$  
of pairwise incomparable prime ideals of $R$ is independent of finite  
character (Theorem \ref{stsubint}). 
 
Conversely, assume that $R$ is a strongly divisorial domain satisfying the  
ascending chain condition on  
prime ideals and that each family $\La$  
of pairwise incomparable prime ideals of $R$ is independent of finite  
character. 
 
Then each subintersection $T$ of $R$  
is of type  $R_{\mc F(\La)}$, where $\La$ is a family of pairwise  
incomparable prime ideals independent of finite  
character. Thus $\mc F(\La)$ is finitely generated \cite [Lemma 1.16]{gab2} and 
$T$ is strongly $w$-divisorial and $t$-flat by Theorem 
\ref{subint}. 
We conclude that $T$ is (strongly) divisorial if and only if each maximal 
ideal of $T$ is a $t$-ideal (Proposition \ref{P1}) if and only if $T$ is flat.  
 
We observe that in 
general, if ${\mc F}$ is a finitely generated multiplicative system  
of ideals, then $R_{\mc F}$  
need not be a flat extension of $R$ \cite[pag. 32]{Fo}. On the other  
hand, we do not know any  
example of a strongly divisorial domain $R$ with a finitely generated 
multiplicative system ${\mc F}$ such that $R_{\mc F}$ is not flat. 
\end{remark}

If $R$ is any domain, we say that $\spec(R)$ (resp.,  
$\tspec(R)$) is {\it treed} (under inclusion) if any maximal (resp.,  
$t$-maximal) ideal of $R$  
cannot contain two incomparable primes (resp., $t$-primes). 
The  
Spectrum of a  Pr\"ufer  domain and the $t$-Spectrum of a $PvMD$ are  
treed.  
If $\spec(R)$ is treed,  
then $\spec(R)= \tspec(R)$ \cite[Proposition 2.6]{KP};  
in particular each maximal ideal is a $t$-ideal and so 
$w$-divisoriality coincides with  
divisoriality by Proposition \ref{P1}.  
 
If $\tspec(R)$ is treed and $\tmax(R)$ is independent of finite  
character,  
then each  
family $\La$ of pairwise incomparable $t$-prime ideals 
of $R$ is independent of finite character. 
 Hence the next  results are easy consequences of  
Theorem \ref{subint} and Theorem \ref{stsubint} respectively. 
 
\begin{corollary} 
     \label{treed} Let $R$ be an integral domain such that  
     $t\op{-Spec}(R)$ is  
treed. The following conditions are equivalent: 
\begin{itemize} \item[(1)] 
$R$ is strongly $w$-divisorial;  
\item[(2)] $R_{{\mc F}(\La)}$ is a $t$-flat $w$-divisorial domain,  
for each set $\La$ of pairwise incomparable $t$-primes; 
\item[(3)] $R_{{\mc F}(\La)}$ is  a $t$-flat strongly $w$-divisorial  
domain,  
for each set $\La$ of pairwise incomparable $t$-primes. 
\end{itemize} 
\end{corollary} 
 
If $R$ has $t$-dimension one, then clearly $\tspec(R)$ is treed. In this  
case, The conditions stated in Corollary \ref{treed} are all satisfied 
if $R$ is $w$-divisorial (cf. Theorem \ref{tmaxtsub}).

\begin{corollary}\label{treed2} Let $R$ be an integral domain such that  
     $t\op{-Spec}(R)$ is  
treed. The following conditions are equivalent: 
\begin{itemize} 
 \item[(1)] $R$ is a strongly $w$-divisorial domain which satisfies the ascending chain conditions on $t$-prime ideals;  
 \item[(2)] Each $t$-subintersection of $R$ is $t$-flat and strongly $w$-divisorial.  
\end{itemize} 
\end{corollary}

%%%%%%%%%%%%%%%%%%%%%%%%%% INTEGRALLY CLOSED DOMAINS  
%%%%%%%%%%%%%%%%%%%%%%%%%%%%% 

\section{Integrally closed $w$-divisorial domains} 
 
 W. Heinzer proved in \cite{H} that an  
integrally closed domain is divisorial if and only if it is 
an $h$-local  Pr\"ufer  domain  
with invertible maximal ideals. We start this section by showing that 
integrally closed $w$-divisorial  
domains have a similar characterization among $PvMD$s.  
Note that a divisorial $PvMD$ is a  Pr\"ufer  domain. 
 
  \begin{lemma} Let $R$ be a $w$-divisorial domain and  $M\in  
t\op{-Max}(R)$. The following conditions are equivalent:  
\begin{itemize}      
\item[(1)] $M$ is $t$-invertible;  
\item[(2)] $MR_M$ is a principal ideal;  
\item[(3)] $R_M$ is a valuation domain.  \end{itemize}  
\end{lemma}  
 
 \begin{proof} $(1)\lra(2)$. Since $t\op{-Max}(R)$ has $t$-finite  
character (Theorem \ref{Loc}), we can apply  
\cite[Theorem 2.2 and Proposition 3.1]{Z}. 
 
  $(2)\ra(3)$  follows from \cite[Lemme 1, Section 4]{Q}, because 
 $R_M$ is a divisorial domain (Theorem \ref{Loc}), and $(3)\ra(2)$  
follows from  
\cite[Lemma 5.2]{H}. \end{proof}  
 
\begin{proposition}\label{corval} Let  
$R$ be a $w$-divisorial domain. Then $R$ is a $PvMD$ 
if and only if each $t$-maximal ideal of $R$  
is $t$-invertible.      
\end{proposition}   
 
\begin{theorem}\label{ic} Let $R$ be an integral  
domain. The following conditions are equivalent:  
\begin{itemize}      
\item[(1)] $R$ is an integrally closed $w$-divisorial domain;  
\item[(2)]  $R$ is a weakly Matlis $PvMD$ and each $t$-maximal ideal  
of $R$ is $t$-invertible.  
\end{itemize}  
\end{theorem} 
 
\begin{proof} $(1)\ra(2)$. A domain $R$ is a $PvMD$ if and  
only if $R$ is an integrally closed $TW$-domain \cite[Theorem 3.5]{K}.  
Hence an integrally closed  $w$-divisorial domain is a $PvMD$.  By  
Theorem \ref{Loc}, $R$ is  a weakly Matlis domain and by Proposition  
\ref{corval} each $t$-maximal ideal is $t$-invertible.  
 
$(2)\ra(1)$. A $t$-maximal ideal $M$ of a $PvMD$ is $t$-invertible  if  
and only if $MR_M$ is a principal ideal \cite{Ho}. Since $R_M$ is  a  
valuation domain, this means that $R_M$ is divisorial \cite[Lemma  
5.2]{H}. Now we can apply Theorem \ref{Loc}. \end{proof} 
 
 The previous theorem can be  proved also by using  the fact that a  
domain $R$ is a $PvMD$ if and only if $R$ is an integrally closed  
$TW$-domain \cite[Theorem 3.5]{K} and the characterization of $PvMD$s  
which are  $TV$-domains given in \cite[Theorem 3.1]{HZ}.  
 
Recall that a  Pr\"ufer  domain $R$ is  strongly discrete if $P^2\neq P$ 
for each nonzero prime ideal $P$ of $R$ \cite[Section 5.3]{fhp} and that  
a generalized Dedekind domain is a strongly discrete  Pr\"ufer  
domain with the property that each ideal has finitely many minimal  
primes \cite{P}.  
We say that a $PvMD$  $R$ is {\it strongly discrete} if  
$(P^2)_t\not=P$, for each $P\in\tspec(R)$ \cite[Remark 3.10]{elb1}.  
If $R$ is a strongly discrete $PvMD$ 
and each $t$-ideal of $R$ has only finitely many minimal primes,  
then $R$ is called a {\it generalized Krull domain} \cite{elb1}.  
 
The next theorem shows 
that the class of strongly $w$-divisorial domains and the class of  
strongly discrete $PvMD$s are strictly related to each other.  
 
\begin{lemma} \label{sd} Let $R$ be a domain. The following conditions 
     are equivalent: 
     \begin{itemize} 
         \item[(1)] $R$ is a strongly discrete $PvMD$; 
         \item[(2)] $R_{M}€$ is a strongly discrete valuation domain, for each  
$M\in t\op{-Max}(R)$; 
\item[(3)] $R_{P}€$ is a strongly discrete valuation domain, for each $P\in 
t\op{-Spec}(R)$; 
\item[(4)] $R_{P}€$ is a valuation domain and $PR_{P}€$ is a principal  
ideal, for each $P\in t\op{-Spec}(R)$; 
\item[(5)] $R_{P}€$ is a divisorial valuation domain,  
for each $P\in t\op{-Spec}(R)$. 
\end{itemize} 
\end{lemma}  
 
\begin{proof} $(1)\lra(4)$. For each $t$-prime ideal $P$ of $R$, we have 
     $(P^2)_{t}= P^2R_{P}€\cap R$ \cite[Proposition  
1.3]{Ho}. Hence $(P^2)_t\neq P$ if and only if  $P^2R_P\neq PR_{P}$.  
Now recall that a maximal ideal of a valuation domain is not idempotent  
if and only if it is principal. 
 
 $(2)\lra(3)$ because each overring of a strongly discrete valuation domain is 
a strongly discrete valuation domain \cite[Proposition 5.3.1(3)]{fhp}.

     $(3)\lra(4)$ by \cite[Proposition 5.3.8 $((2) \lra (6))$]{fhp}. 
    
$(4)\lra(5)$ by \cite[Lemma 5.2]{H}. 
     \end{proof}

\begin{comment} 
\begin{proposition}\label{ics} Let $R$ be an integral domain. The  
following conditions are equivalent: 
 
\begin{itemize} 
 
\item[(1)] $R$ is an integrally closed strongly $w$-divisorial domain; 
 
\item[(2)] $R$ is a strongly discrete $PvMD$ and a weakly Matlis  
domain. 
 
\end{itemize} \end{proposition} 
 
\begin{proof} $(1)\ra(2)$. By Theorem \ref{ic}, $R$ is a weakly  
Matlis $PvMD$. Now, given $P\in\tspec(R)$, $R_P$ is a divisorial  
valuation domain. Hence 
 $R$ is a strongly discrete $PvMD$ by Lemma \ref{sd}. 
 
$(2)\ra(1)$. Clearly $R$ is integrally closed. In addition by  
Lemma \ref{sd} 
$R_{P}€$ is a divisorial domain, for each $P\in\tspec(R)$ . 
Hence $R$ is a strongly $w$-divisorial domain. 
\end{proof} 
\end{comment}

\begin{theorem}\label{twic} Let $R$ be an integral domain. The  
following conditions are equivalent:  
 
\begin{itemize}  
 
\item[(1)] $R$ is a strongly discrete $PvMD$ and a weakly Matlis  
domain; 
 
\item[(2)] $R$ is an integrally closed strongly $w$-divisorial domain; 
 
\item[(3)] $R$ is  integrally closed and each $t$-flat overring of $R$  
is  $w$-divisorial; 
 
\item[(4)] $R$ is integrally closed and each $t$-linked overring of $R$  
is $w$-divisorial; 
 
\item[(5)] $R$ is a $w$-divisorial generalized Krull domain; 
 
\item[(6)] $R$ is a generalized Krull domain and each $t$-prime  
ideal of $R$ is contained in a unique $t$-maximal ideal. 
 
\end{itemize} \end{theorem} 
 
 \begin{proof} 
     $(1)\ra(2)$. Clearly $R$ is integrally closed. In addition, by  
Lemma \ref{sd}, 
$R_{P}€$ is a divisorial domain, for each $P\in\tspec(R)$ . 
Hence $R$ is a strongly $w$-divisorial domain. 
 
$(2)\ra(3)$. By Theorem \ref{ic}, $R$ is a $PvMD$; in particular $\tspec(R)$ 
is treed. Thus we can apply Corollary  \ref{treed}. 
 
$(3)\ra(1)$. By Theorem \ref{ic}, $R$ is a weakly  
Matlis $PvMD$. Now, given $P\in\tspec(R)$, $R_P$ is a divisorial  
valuation domain. Hence 
 $R$ is a strongly discrete $PvMD$ by Lemma \ref{sd}. 
 
$(3)\lra(4)$. By Theorem \ref{ic}, statements  
(3) and (4) imply that $R$ is a $PvMD$. The conclusion now follows from  
the  
fact that each $t$-linked overring of a $PvMD$ $R$ is $t$-flat  
\cite[Proposition 2.10]{KP}. 
 
$(1)\ra(5)$. By (1)$\ra$(2), $R$ is a  
$w$-divisorial domain. To show that $R$ is a generalized Krull  
domain, let $I$ be a  
$t$-ideal of $R$. Since $R$ has $t$-finite character, then $I$ is  
contained in only finitely many $t$-maximal ideals. Furthermore, each  
$t$-prime ideal is contained in a unique $t$-maximal ideal. Thus $I$  
has  
just finitely many minimal ($t$)-prime ideals. We conclude by  
using \cite[Theorem 3.9]{elb1}. 
 
$(5)\ra(6)$ is clear. 
 
$(6)\ra(1)$. It is enough to show that $R$ has $t$-finite character. 
This follows from the fact that each nonzero principal ideal has 
finitely many minimal ($t$)-primes. 
\end{proof} 
 
As a consequence of Theorem \ref{twic}, we obtain the following  
characterization of  integrally closed totally divisorial domains (see also  
\cite{O}). 
 
\begin{corollary}\label{tvic} Let $R$ be an integral domain. The  
following conditions are equivalent:  
 
\begin{itemize}  
 
\item[(1)] $R$ is an integrally closed  totally divisorial domain; 
 
\item[(2)] $R$ is  integrally closed and each flat overring of $R$ is  
divisorial; 
 
\item[(3)] $R$ is an  integrally closed strongly divisorial domain; 
 
 \item[(4)] $R$ is an $h$-local strongly discrete  Pr\"ufer  domain; 
 
\item[(5)] $R$ is a divisorial generalized Dedekind domain; 
 
\item[(6)] $R$ is a generalized Dedekind domain and each nonzero prime  
ideal is contained in a unique maximal ideal. 
 
\end{itemize} \end{corollary} 
 
\begin{proof} This follows from the fact that in a  Pr\"ufer  domain the  
$d$- and $t$-operation coincide, that each overring of a  Pr\"ufer  
domain is a flat  Pr\"ufer  domain, and that a  Pr\"ufer  domain is  
a generalized Krull domain if and only if it is a 
generalized Dedekind domain \cite{elb1}. \end{proof}

Recall that the {\it complete integral closure} of $R$ is the overring  
$\widetilde{R} := \cup\{(I\colon I)\; ; \; I$ nonzero ideal of $R\}$.  
If $R = \widetilde{R}$, we say that $R$ is {\it completely integrally  
closed}. 
 
 \begin{proposition}\label{P2} Let $R$ be an integral domain. The  
following conditions are equivalent: 
     \begin{itemize} 
     \item[(1)] $R$ is an integrally closed $w$-divisorial domain of  
$t$-dimension one; 
  \item[(2)] $R$ is an integrally closed domain of  
$t$-dimension one and each $t$-linked overring of $R$ is  
$w$-divisorial; 
  \item[(3)] $R$ is a completely integrally closed $w$-divisorial  
domain; 
  \item[(4)] $R$ is a Krull domain. 
\end{itemize} \end{proposition} 
\begin{proof} $(1)\lra(2)\lra(4)$. Clearly a $w$-divisorial domain of  
$t$-dimension one is strongly $w$-divisorial. Since  
a generalized Krull domain of  
$t$-dimension one is a  Krull domain \cite[Theorem  
3.11]{elb1}, we can conclude by  
applying Theorem \ref{twic}. 
 
$(3)\lra(4)$ because a  
completely integrally closed $TV$-domain is Krull \cite[Theorem  
2.3]{HZ}.     
\end{proof} 
 
  It is well-known that a divisorial Krull domain is a Dedekind domain;  
hence by the previous proposition  we recover that a completely  
integrally closed  divisorial domain is a Dedekind domain  
\cite[Proposition 5.5]{H}.   
 
\begin{remark} \label{cic} \rm  
     Recall that, for any domain $R$, $\widetilde{R}$ is  
integrally closed and $t$-linked over $R$ \cite[Corollary 2.3]{DHLZ}.  
Since each localization of a $t$-linked overring of $R$  
is still $t$-linked over $R$, if each  
$t$-linked overring of $R$ is $w$-divisorial, we have that 
$\widetilde{R}$ is an integrally closed strongly $w$-divisorial domain. 
In this case, by Theorem \ref{twic}, $\widetilde{R}$  is  
a weakly Matlis strongly discrete $PvMD$. 
If in addition $\widetilde{R}$ is completely integrally closed, for  
example 
if $(R\colon \widetilde{R}) \neq 0$, by Proposition \ref{P2}  
$\widetilde{R}$ is  
a Krull domain. 
 
In a similar way, by using Corollary \ref{tvic}, we see that if $R$ is  
totally divisorial, the integral closure of $R$ is an $h$-local  
strongly divisorial  Pr\"ufer  domain. 
\end{remark}

%%%%%%%%%%%%%%%%%%%%%%%%%%%%%%%%%% MORI DOMAINS  
%%%%%%%%%%%%%%%%%%%%%%%%%%%%%%%%%%% 

\section{Mori $w$-divisorial domains}

We start by recalling some properties of Noetherian divisorial  
 domains proved in \cite{H, Q}.  
 
\begin{proposition} \label{Noeth1} Let $R$ be a  
 domain. The following conditions are equivalent: 
\begin{itemize} \item[(1)] $R$ is a one-dimensional $w$-divisorial Mori 
domain;  
\item[(2)] $R$ is a divisorial Mori domain; 
\item[(3)]  $R$ is a divisorial Noetherian domain; 
\item[(4)] $R$ is a Mori domain and each two generated  
ideal of $R$ is divisorial; 
\item[(5)] $R$ is a one-dimensional Mori domain and $(R\colon  
     M)$ is a two generated ideal, for each $M\in \op{Max}(R)$; 
\item[(6)] $R$ is a one-dimensional Noetherian domain and $(R\colon  
     M)$ is a two generated ideal, for each $M\in \op{Max}(R)$. 
\end{itemize}  
     \end{proposition}  
     \begin{proof} $(1)\ra(2)$ by Proposition \ref{P1}. 
 
$(2) \ra(3)$ because each $v$-ideal of a Mori domain  
     is $v$-finite. 
 
$(3)\ra(1)$ because 
     Noetherian divisorial domains are one-dimensional \cite[Corollary  
4.3]{H}. 
 
     $(3) \lra (6)$ and $(2)\lra(4)\lra(5)$ by \cite[Theorem 3, Section 2]{Q}. 
\end{proof} 
 
An integrally closed $w$-divisorial Mori domain is 
a Krull domain. In fact it has to be a $PvMD$ (Theorem \ref{ic}).  
By Proposition \ref{Noeth1}, any Noetherian integrally closed  
domain of dimension greater 
than one is a $w$-divisorial Noetherian domain that is not divisorial.  
 
We say that a nonzero fractional ideal $I$ of $R$ is a {\it $w$-divisorial  
ideal} if $I_v=I_w$. With this notation, a $w$-divisorial domain  
is a domain in which each nonzero ideal is $w$-divisorial.  
We also say that, for $n \geq 1$,  $I$ is $n$ $w$-generated 
if $I_{w}€=(a_{1}R+ \dots+a_{n}R)_{w}€$, for some $a_{1}, \dots,a_{n}$  
in the quotient field of $R$.

\begin{theorem} \label{wMori} Let $R$ be a Mori  domain. The  
following conditions are equivalent: 
 
 \begin{itemize} \item[(1)] $R$ is a $w$-divisorial domain; 
 
\item[(2)] Each two generated  nonzero ideal is $w$-divisorial; 
 
\item[(3)] $R$ has $t$-dimension one and $(R\colon M)$ is a two  
$w$-generated  
ideal, for each $M \in t\op{-Max}(R)$.  \end{itemize} 
 
\end{theorem} 
 
\begin{proof}$(1)\ra(2)$ is clear. 
 
$(2)\ra(3)$. Let $M\in \tmax(R)$. Since $R$ is a Mori domain, then $M$  
is a divisorial ideal. Let $x\in (R\colon M)\sm R$, then $(R\colon  
M)=(R+Rx)_v$.  
So that by assumption $(R\colon M)=(R+Rx)_w$. To conclude, we show  
that  
$R_M$ is one-dimensional. Let $I$ be a nonzero two generated ideal of  
$R_M$. 
Then,  we can assume that $I=(a, b)R_M$ for some $a, b\in I\cap R$.  
Since $R$  
is a Mori domain, then $I_v=((a, b)R_M)_v=(a, b)_vR_M$. Hence  
$I_v=(a,b)_wR_M=(a, b)R_M=I$.  
Thus each two generated ideal of $R_M$ is  
divisorial. It follows from Proposition \ref{Noeth1} that  
$R_M$ is one-dimensional. 
 
$(3)\ra(1)$.  Since $R$ is a $TV$-domain, by Theorem \ref{Loc},  
it is enough to  
show that $R_{M}€$ is a divisorial domain for each $M\in\tmax(R)$. 
This follows again from Proposition \ref{Noeth1}. In fact, 
by assumption $R_M$ is a Mori domain of dimension one.  
Let $(R\colon M)=(a, b)_w$ for some $a, b\in (R\colon M)$. Then  
$(R_{M}€\colon MR_{M}€)=(R\colon M)R_M=(a, b)_wR_M=(a, b)R_M$  is  
two generated  
(the first equality holds because $M$ is $v$-finite). 
\end{proof} 
 
Recall that a {\it Strong Mori domain} is a domain satisfying the  
ascending chain condition on $w$-ideals.  
A domain $R$ is a Strong Mori domain if and only if it has $t$-finite  
character and $R_{M}€$ is Noetherian,  
for each $t$-maximal ideal $M$ \cite[Theorem 1.9]{FM2}. Thus a Mori  
domain is Strong Mori if and only if $R_{M}€$ is Noetherian,  
for each $t$-maximal ideal $M$. 
 
\begin{corollary}\label{SM} \cite[Corollary 2.5]{M2}     
     A $w$-divisorial Mori domain 
   is a Strong Mori domain of $t$-dimension one. 
\end{corollary}  
\begin{proof}  
     A $w$-divisorial Mori domain is Strong Mori (because $w=v$) 
         and has $t$-dimension one by Theorem \ref{wMori}.  
     \end{proof} 
     
We next investigate $w$-divisoriality of overrings of Mori  
domains. Our first result in this direction shows that, if $R$ is Mori,  
$w$-divisoriality is inherited by generalized ring of fractions.  
This improves \cite[Theorem 2.4]{M2}. 
 
We observe that a Mori domain is a $v$-coherent $TV$-domain,  
because each $t$-ideal of a Mori domain is $v$-finite. 
 We also recall that if $R$ is 
$v$-coherent, we have $I_tR_S=(IR_S)_t$ , for each nonzero fractional ideal $I$ 
and  each multiplicative set $S$. 
 
 \begin{proposition}\label{TW} Let $R$ be a $v$-coherent domain. The  
following conditions 
are equivalent: 
 
\begin{itemize} 
 
     \item[(1)] $R$ is a $TW$-domain; 
 
     \item[(2)]  All the nonzero ideals of $R_{M}$ are $t$-ideals, for each $M\in  
     t$-$\op{Max}(R)$; 
 
     \item[(3)]  All the nonzero ideals of $R_{P}€$ are $t$-ideals, for each $P\in  
     t$-$\op{Spec}(R)$; 
     
\item[(4)] Each $t$-flat overring of $R$ is a $TW$-domain. 
 
\end{itemize} 
\end{proposition} 
 
\begin{proof}  
 
     $(1)\lra(2)$. Let $I$ be a nonzero ideal and $M$ a $t$-maximal ideal of $R$.  
     If $t = w$ on $R$, then 
 $IR_{M}=I_{w}R_{M}=I_{t}R_{M}=(IR_{M})_{t}$. 
 
 Conversely, we have 
$IR_{M}€=(IR_{M})_{t}=I_{t}€R_{M}$. 
Thus  
$$I_w=\cap_{M\in\tmax(R)}IR_M=\cap_{M\in\tmax(R)}I_tR_M=I_t.$$  
 
$(2)\ra(3)$. Let $I$ be a nonzero ideal of $R$, $P$ a $t$-prime of $R$ and $M$  
a  
$t$-maximal ideal containing $P$. Then 
$$IR_{P}=(IR_M)R_P= 
(IR_M)_tR_P=(I_{t}R_{M})R_P=I_{t}R_{P}=(IR_{P})_{t}.$$ 
 
$(3)\ra(4)$. Let $T$ be a $t$-flat overring of $R$. Then 
 $T$ is a $v$-coherent domain \cite[Proposition 3.1]{gab2}.  
If $N$ is a $t$-maximal ideal of $T$, then $P = N\cap R$ 
is a $t$-prime of $R$ and $T_{N}= R_{P}$. Hence, 
 if (3) holds, each nonzero ideal of $T_{N}$ is  
a $t$-ideal and $T$ is a $TW$-domain by  $(2)\ra(1)$.  
 
$(4)\ra(1)$ is clear.  
\end{proof} 
 
\begin{theorem}\label{wvMori} Let $R$ be a Mori domain. 
     The following conditions are equivalent: 
 
 \begin{itemize}      
 
\item[(1)] $R$ is  $w$-divisorial;      
 
\item[(2)] $R$ is  strongly $w$-divisorial; 
 
\item[(3)] Each $t$-flat overring of $R$ is $w$-divisorial;  
 
\item[(4)]  Each generalized ring of fractions of $R$ is  
$w$-divisorial;   
 
\item[(5)] $R_{M}$ is a divisorial domain, for each  $M\in  
t\op{-Max}(R)$.  
 
\end{itemize} \end{theorem}   
 
\begin{proof}  
Each generalized ring of fractions of a Mori domain is  
Mori  \cite[Corollaire 1, Section 3]{Q}; thus it is a $TV$-domain.  
In addition,  
each generalized ring of fractions of a Mori domain is  
$t$-flat, because each $t$-ideal is $v$-finite and so each multiplicative 
system of ideals is $v$-finite. Hence we can apply Proposition \ref{TW}. 
     \end{proof} 
 
$t$-linked overrings of Mori domains do not behave as well as  
generalized rings of fractions. In fact a Mori non-Krull  
domain has $t$-linked overrings which are not 
$t$-flat \cite[Corollary 2.10]{DHLZ2}.  
Also, if each  
$t$-linked overring of a Mori domain $R$ is Mori, then $R$ has  
$t$-dimension one \cite[Proposition 2.20]{DHLZ}. The converse holds  
if $R$ is a Strong Mori domain; precisely, we  
have the following result. 
 
\begin{proposition}\label{overSM} 
     Each $t$-linked overring of a Strong Mori domain  
of $t$-dimension one is either a field or a Strong Mori domain  
of $t$-dimension one.  
\end{proposition} 
\begin{proof}  
It follows from \cite[Theorem 3.4]{FM2} recalling that  
an overring of a domain is a $w$-module if and only if it is  
$t$-linked \cite [Proposition 2.13 (a)]{DHLZ}. 
\end{proof}

\begin{corollary}\label{coroverSM} If $R$ is a $w$-divisorial Mori  
domain, then each 
$t$-linked overring of $R$ is either a field or a Strong Mori domain of $t$-dimension  
one.  
\end{corollary} 
\begin{proof} 
It follows from Corollary \ref{SM} and Proposition \ref{overSM}. 
\end{proof} 
 
Our next purpose is to improve and generalize to Mori domains some  
results  
proved in \cite{BS} for Noetherian totally divisorial domains. 
 
\begin{proposition}\label{Noeth2} Let $R$ be a domain.  
     The following conditions are equivalent: 
\begin{itemize} 
\item[(1)] $R$ is a one-dimensional domain and each 
$t$-linked overring of $R$ is $w$-divisorial; 
\item[(2)] $R$ is a one-dimensional totally divisorial domain; 
 \item[(3)] $R$ is a Noetherian totally divisorial domain; 
 \item[(4)] Each ideal of $R$ is two generated. 
\end{itemize} 
\end{proposition} 
\begin{proof} $(1)\ra(2)$. Since dim$(R) =1$, each overring  
of $R$ is $t$-linked over $R$ \cite[Corollary 2.7 (b)]{DHLZ}. 
Hence each overring $T$ of $R$ is  
$w$-divisorial. Assume that $T$ is not a field. 
To prove that $T$ is divisorial it suffices to check  
that dim$(T)=1$ (Proposition \ref{P1}). Let $R'$ be the integral  
closure  
of $R$ and $T'$ that of $T$.  
Since $R'$ is one-dimensional and $w$-divisorial, then 
$R'$ is divisorial. Thus $R'$, being integrally closed, is a  Pr\"ufer  
domain  
\cite[Theorem 5.1]{H}. It follows that the extension  
$R'\subseteq T'$ is flat, and so dim$(T')\le$ dim$(R')=1$. Hence 
dim$(T)=$ dim$(T')=1$. We conclude that  $T$ is divisorial and  
therefore $R$ is  
totally divisorial. 
     
     $(2)\ra(3)$ by \cite[Proposition 7.1]{BS}. 
     
     $(3)\ra(1)$ by Proposition \ref{Noeth1}. 
     
     $(3)\lra(4)$ by \cite[Theorem 7.3]{BS}, because in the Noetherian case a 
     domain is totally divisorial if and only if it is totally reflexive \cite[Section 3]{O2}. 
     \end{proof}

Lemma \ref{ngen} below is similar to \cite[Theorem 26(2)]{Ma}. 
We will need the following version  
of Chinese Remainder Theorem, whose proof is straightforward. 
 
\begin{lemma} \label{CRT}  Let $R$ be an integral domain, $I$ an ideal  
of $R$,  
$P_1,\ldots,  P_n$  a set of pairwise incomparable prime ideals and $S =  
R\sm (P_{1}€\cup\dots \cup P_{n})€$. 
If $x_1,\ldots, x_n\in I$, there exists $x\in IR_{S}€$ such that 
$x\equiv x_i\, (\op{mod}\, IP_iR_{P_i})$, for each $i=1,\ldots, n$. 
\end{lemma} 
 
\begin{lemma}\label{ngen} Let $R$ be an integral domain which has  
$t$-finite character and $I$ a nonzero ideal of $R$. Let $n$ be a  
positive integer 
and assume that, for each 
$M \in t\op{-Max}(R)$, a minimal set of generators of 
$IR_M$ has at most $n$ elements. Then $I$  
is $w$-generated by a number of generators $m \le\max(2,n)$. 
\end{lemma} 
 
\begin{proof} If $I$ is not contained in any $t$-maximal ideal, then  
$I_{w}€= R$. Otherwise, let $M_1,\ldots, M_r$  
be the $t$-maximal ideals of $R$ which contain $I$. For $i=1,\ldots,  
r$, let $a_{1i},\ldots,a_{ni}\in I$ be such that  
$IR_{M_i}=(a_{1i},\ldots,a_{ni})R_{M_i}$.  By Lemma \ref{CRT}, if 
$S =  
R\sm (M_{1}€\cup\dots \cup M_{r})€$, for each  
$j=1,\ldots, n$, there exists $a_j\in IR_{S}€\sub IR_{M_{i}€}€$ such 
that $a_j\equiv a_{ji}(\op{mod} \, IM_iR_{M_i})$, for each $i=1,\ldots, r$. 
By going modulo $IM_iR_{M_i}$ and using Nakayama's Lemma, we  
get $IR_{M_i}=(a_1,\ldots, a_n)R_{M_i}$ for each $i=1,\ldots, r$. We  
can assume that the $a_{j}$'s are in $I$ and $a_1\not=0$. Let  $N_1, \ldots, N_s$ 
be the set of  
$t$-maximal ideals which contain $a_1$, with $N_1=M_1, \ldots, N_r=M_r$. Let  
$b\in I\setminus\cup_{j=r+1}^sM_j$. Then $IR_{N_j}=(a_1,\ldots,  
a_n)R_{N_j}$ for $j=1,\ldots, r$ and $IR_{N_j}=(a_1, b)R_{N_j}  
=R_{N_j}$ for $j=r+1,\ldots, s$. By arguing as above, there exist  
$b_1=a_1, b_2, \ldots, b_n\in I$ such that $IR_{N_j}=(b_1,\ldots,  
b_n)R_{N_j}$ for each $j=1,\ldots, s$. We claim that  $I_w=(b_1,\ldots,  
b_n)_w$. Let M be a $t$-maximal ideal of $R$. If $M=N_j$ for some $j$,  
then  $IR_M=(b_1,\ldots, b_n)R_M$. If $M\not=N_j$ for $j=1,\ldots, s$,  
then  $IR_M=R_M=(b_1,\ldots, b_n)R_M$, since $b_1=a_1\notin M$. 
     \end{proof}

\begin{theorem}\label{Mori2} Let $R$ be a domain. The  
following conditions are equivalent: 
 
\begin{itemize} 
 
\item[(1)] $R$ has $t$-dimension one and  
each $t$-linked overring of $R$ is $w$-divisorial; 
 
\item[(2)] $R$ is a Mori domain and  
each $t$-linked overring of $R$ is $w$-divisorial; 
 
\item[(3)] $R$ is a Mori domain and $R_M$ is totally divisorial, for  
each  
$M \in t\op{-Max}(R)$; 
 
 \item[(4)]  Each nonzero ideal of $R$ is a two $w$-generated $w$-divisorial ideal;
 
\item[(5)]  Each nonzero ideal of $R$ is two $w$-generated. 
 
\end{itemize} 
\end{theorem} 
 
\begin{proof} $(1)\ra(2)$.  
     $R$ has $t$-finite character, because it is $w$-divisorial (Theorem  
     \ref{Loc}). 
We now show that, for each $M\in\tmax(R)$, $R_M$ is Noetherian.  
Since $R_{M}€$ is a one-dimensional  
$t$-linked overring of $R$, then $R_{M}€$ is divisorial (Proposition  
\ref{P1}). In addition, each overring $T$ of $R_{M}€$ is $t$-linked  
over $R_{M}€$ \cite[Corollary 2.7]{DHLZ} and so it is $t$-linked over  
$R$. Thus $T$ is a $w$-divisorial domain. By Proposition \ref{Noeth2},  
$R_M$ is Noetherian. We conclude that $R$ is a (Strong) Mori domain. 
     
$(2)\ra(3)$.     
     $R$ is clearly $w$-divisorial. Hence $R_M$  
is a one-dimensional Noetherian domain (Corollary \ref{SM}). Let $T$ be 
a $t$-linked overring of $R_M$. Hence  
$T$ is $t$-linked over $R$ and so by assumption it is  
$w$-divisorial. By Proposition \ref{Noeth2} 
$R_M$ is totally divisorial. 
 
$(3)\ra(4)$. $R$ is $w$-divisorial by Theorem \ref{wvMori}. 
Hence $R_{M}$ is  one-dimensional and Noetherian by Corollary \ref{SM}. 
Thus, for each $M\in \tmax(R)$, each ideal of  
$R_M$ is two generated by Proposition \ref{Noeth2}.  
By using Lemma \ref{ngen}, we 
conclude that  every nonzero ideal of $R$ is  a two $w$-generated $w$-divisorial ideal. 
 
 $(4)\ra(5)$ is clear.
 
$(5)\ra(3)$. If $(5)$ holds, 
$R$ is a Strong Mori domain and so $R_M$ is a Noetherian domain, for each  
$M\in\tmax(R)$. Let $IR_M$ be a nonzero ideal of $R_M$, where $I$ is an  
ideal of $R$. By assumption, $I_w=(a,b)_w$ for some $a, b \in R$.  
Thus $IR_M = (a,b)_wR_M =(a,b)R_M$ is a two generated ideal.  
 It follows from Proposition \ref{Noeth2} that $R_M$ is  
 a totally divisorial domain.  
 
 $(3)\ra(2)$. $R$ is $w$-divisorial by Theorem \ref{wvMori}. 
Let $T$ be a  
$t$-linked overring of $R$, $T \neq K$. By Corollary \ref{coroverSM}, 
$T$ is a  Mori domain. To show that $T$ is $w$-divisorial, 
by Theorem \ref{wvMori}, we 
have to prove that $T_N$ is a divisorial domain, for each $N\in  
\tmax(T)$. 
 Since $R\sub T$ is $t$-linked, then $Q =(N\cap R)_t\not=R$  
\cite[Proposition 2.1]{DHLZ}; but as $R$ has $t$-dimension one  
(Corollary \ref{SM}), then $Q$ is a  
$t$-maximal ideal of $R$. Since $R_Q$ is totally divisorial and 
$R_Q \sub T_N$, then $T_N$ is a divisorial domain.  
 
$(2)\ra(1)$ by Corollary \ref{SM}. 
\end{proof} 
 
 \begin{corollary}  Let $R$ be a domain and assume that each $t$-linked overring of $R$ is $w$-divisorial. Then $R$ is a Mori domain if and only if it has $t$-dimension one.
 \end{corollary}

 \begin{example} \rm  
Mori non-Krull and non-Noetherian domains satisfying the equivalent 
conditions of Theorem \ref{Mori2} can be constructed by using pullbacks, as the  
following example shows.

Let $T$ be a Krull domain having a maximal ideal $M$ of height one and 
assume that the residue field $K =T/M$ has a subfield $k$ such that  
$[K:k]=2$. Let $R=\varphi^{-1}(k)$ be the pullback of $k$ with  
respect 
to the canonical projection $\varphi : T\longra K$.  
 
The domain $R$ is Mori and it is Noetherian  
if and only if $T$ is Noetherian \cite[Theorems 4.12 and 4.18]{GH1}.  
$M$ is a maximal ideal of $R$ that is divisorial;  
thus $M\in \tmax (R)$. Since $R_{M}€$ is the  
pullback of $k$ with respect to the natural projection $T_{M}€\longra  
K$,  
$R_{M}€$ is divisorial by \cite[Corollary 3.5]{M2}. 
In addition $T_{M}€$ is the only overring of $R_{M}€$. In fact  
each overring of $R_{M}€$ is comparable with $T_{M}€$ under inclusion;  
but $T_{M}€$ is a $DVR$ and $[K:k]=2$.  
Thus $R_{M}€$ is totally divisorial. 
 
If $N$ is a  
$t$-maximal ideal of $R$ and $N\neq M$, there is a unique $t$-maximal  
ideal $N'$ of $T$ such that $N'\cap R=N$ \cite[Theorem 2.6(1)]{GH} and  
for this prime $T_{N'}€=R_{N}€$. Thus $R_{N}€$ is a $DVR$.  
It follows that $R_{N}€$ is totally divisorial, for each $N \in  
\tmax(R)$. 
\end{example}

 \end{document}